\documentclass[11pt]{amsart}

\usepackage{amsmath}
\usepackage{amssymb}
\usepackage{graphicx}

\newtheorem{theorem}{Theorem}[section]

\newtheorem{lemma}[theorem]{Lemma}
\newtheorem{propo}[theorem]{Proposition}
\theoremstyle{definition}

\newtheorem{example}[theorem]{Example}

 \newtheorem{assumption}{Assumption}[section]

\numberwithin{equation}{section}

\usepackage{xpatch} 
\makeatletter   %
\xpatchcmd{\@thm}{\thm@headpunct{.}}{\thm@headpunct{}}{}{}

\def\E{\mathbb{E}}
\def\P{\mathbb{P}}
\def\R{\mathbb{R}}

\newcommand{\eeq}{\end{equation}}

\newcommand{\eeqa}{\end{eqnarray}}

\newcommand{\baa}{\begin{eqnarray*}}
\newcommand{\eaa}{\end{eqnarray*}}

\newcommand{\M} {\tilde{M}}
\newcommand{\G} {\overline{\Gamma}}
\newcommand{\GN} {\tilde{\Gamma}^N}

\title[Spatial SIR epidemic model with varying infectivity]{Spatial SIR epidemic model with varying infectivity without movement of individuals}

\author[A. Kanga]{Armand Kanga$^1$}
\address{$^1$Laboratory of Applied Mathematics and Computer Science, Universit{\'e} Felix Houphou{\"e}t-Boigny,  Abidjan, Ivory Coast and Aix--Marseille Universit{\'e}, CNRS, I2M, Marseille, France }
\email{kangayaoarmand@gmail.com}
\author[{\'E}.Pardoux]{{\'E}tienne Pardoux$^2$}
\address{$^2$ I2M, Aix--Marseille Universit{\'e}, CNRS, Marseille, France}
\email{etienne.pardoux@univ-amu.fr}



\begin{document}



%
%

\begin{abstract} 
	We study an SIR epidemic model with variable infectivity,  where the individuals are distributed over a compact subset 
	$D$ of $\R^d$. We define empirical measures which describe the evolution of the state (susceptible, infectious, recovered) of the individuals in the various locations, and the total force of infection in the population. In our model, the individuals do not move. We establish a law of large numbers for these measures, as the population size tends to infinity.
\end{abstract}

\maketitle

\section{Introduction}
Epidemic models using ordinary differential equations have been the subject of much research in recent years. Anderson and Britton \cite{e}, Britton and Pardoux \cite{d} have shown that these models are limits, when the population size tends towards infinity, of  stochastic Markovian models. In particular, the Markovian nature of this model implies  that the duration of infection is exponentially distributed, which is unrealistic for most epidemics.\\
As a result, models with non-exponential infection durations have attracted some interest, see  in particular \cite{f} and \cite{g}. Kermack and McKendrick \cite{n} also considered that the infectivity should be a function which varies with the time since infection. The duration of infection is the time taken by this function to vanish out definitively; its law is completely arbitrary. In \cite{j}, the authors have established the law of large numbers for the SIR model with variable infectivity, where the infectivity varies from one individual to another and depends upon the time elapsed since infection. They assume that the infectivity function has a finite number of jumps, and satisfies an assumption of uniform continuity between jumps. In \cite{b}, the same law of large numbers is established under a weaker assumption: infectivity functions have their trajectories in $\mathbb{D}(\mathbb{R}_+,\mathbb{R})$, and are bounded by a constant. However, in the various models studied above, the authors ignore the fact that a population extends over a spatial region. Yet, spatial heterogeneity, habitat connectivity and movement rates play an important role in the evolution of infectious diseases. Both deterministic and stochastic models have been used to understand the importance of the heterogeneity of the density of individuals on the spread of infectious diseases, on the persistence or extinction of an endemic disease, for example \cite{r}, \cite{s}, \cite{t}
and \cite{PPdense}. Some Markovian models in this framework have been studied in \cite{a}. They studied a stochastic SIR compartmental epidemic model for a population which moves on a torus $(\mathbb{T}^2 = \mathbb{R}^2/\mathbb{Z}^2)$ according to Stochastic Differential Equations driven by independent Brownian motions. They define sequences of empirical measures that describe the evolution of the positions of susceptible, infected and recovered individuals.  They establish large-population approximations of these sequences of  measures. In \cite{c}, the authors consider a population distributed in the space $\mathbb{R}^d$ whose individuals are characterized by: a position and an infection state, interact and move in $\mathbb{R}^d.$ An epidemic model combining spatial structure and variable infectivity would be more realistic. This is the focus of our work. As a result, we are considering a population distributed over a compact  
 subset $D$ of $\R^d$; and use the same type of arguments as in \cite{b}. We define sequences of empirical measures which describe the evolution of the positions of susceptible, infected and recovered individuals, and  establish the law of large numbers for these measures. In this paper, we restrict ourselves to the case where the individuals do not move. Note however that infectious individuals can infect susceptible individuals located far away. This means that we can take into account movements of individuals (daily from home to work, or occasionally for vacation and visits to the family), without modeling those movements explicitly.
 
 The same model, with diffusive movement of the individuals will be considered in another publication.
\subsection{Notation}
We denote by
\begin{itemize}
	\item $\mathcal{M}$ denotes the set of finite positive measures  on $D$ which we equip with the topology of weak convergence.
	\item $\mathbb{D}:=\mathbb{D}(\mathbb{R}_+,\mathbb{R}_+)$ denotes the space of c\`adl\`ag functions  defined on $\mathbb{R}_+$  with values in  $\mathbb{R}_+$.
	\item $\mathbb{D}_{\mathcal{M}}:=\mathbb{D}(\mathbb{R}_+,\mathcal{M})$ denotes the space of c\`adl\`ag functions  defined on $\mathbb{R}_+$  with values in $\mathcal{M}$.
	\item For all $ \varphi \in C_b(D)$ and  $\mu \in \mathcal{M}$, $(\mu,\varphi)= \displaystyle \int_{D} \varphi(y) \mu(dy)$.
\end{itemize}
 $c$ and $C$ denote positive constants that can change from line to line.

\section{Model description}
\label{k}
The epidemic model studied here is the SIR model in a spatial framework with variable infectivity; the letters S, I and R represent the different states of an individual ("susceptible","infected" and "recovered" respectively). The SIR model states that a  susceptible individual can become infected, and finally  recovered when he/she recovers from the disease. In our spatial model, an individual is characterized by its state $E\in \left\{ S, I, R \right\} $ and its position $X$, a point in $D$ which is a compact subset of $\R^d$. To simplify the mathematical description, we identify the $S$, $I$ and $R$ states as $0$, $1$ and $2$ respectively. The  space of individuals is therefore $D\times \left\{ 0, 1, 2 \right\} $. We consider a population of fixed size N; and we assume that at time t=0 the population is divided into three subsets: those susceptible , there are $S^N(0)$ of them, those infected, there are $I^N(0)$ of them, and those recovered, there are $R^N(0)$ of them i.e $S^N(0)+I^N(0)+R^N(0)= N$.  We denote by $\{X^{i}, i\in \{ \mathfrak{S}_N,\mathfrak{I}_N,\mathfrak{R}_N\}\}$ the positions of the  individuals at time $t=0$, where $(\mathfrak{S}_N,\mathfrak{I}_N,\mathfrak{R}_N)$ forms a partition of $\{1,\cdots,N\}$ with $Card(\mathfrak{S}_N)=S^N(0)$, $Card(\mathfrak{I}_N)=I^N(0)$ and $Card(\mathfrak{R}_N)=R^N(0)$. Of course the three sets
$\mathfrak{S}_N$, $\mathfrak{I}_N$ and $\mathfrak{R}_N$ depend upon $N$.
Now let us consider $\{\lambda_{-j} , j\geq 1\}$ and $\{\lambda_{j} , j\geq 1\}$ two mutually independent sequences of  i.i.d random elements  of $\mathbb{D}(\mathbb{R}_+,\mathbb{R}_+).$ $\lambda_{-j}(t)$ is the infectivity at time $t$ of the individual $j \in\mathfrak{I}_N$ and $\lambda_{j}(t)$ is the infectivity at time $t$ after its infection of the individual $j\in\mathfrak{S}_N$. We assume that there exists a deterministic constant $\lambda^*>0$ such that $0\leq\lambda_{j}(t)\leq \lambda^*$ a.s, for all $j\in \mathbb{Z}^* \text{ and } t\geq 0$, with the convention: $\forall j \ge1$, $\lambda_{j}(t) = 0 \text{ if } t < 0$ and we shall use the notations $\overline{\lambda}^0(t)= \mathbb{E}(\lambda_{-1}(t))$ and $\overline{\lambda}(t)= \mathbb{E}(\lambda_{1}(t))$. It is natural that an infected individual is more likely to infect a close neighbor than a more distant one. While these different transmission behaviors are averaged in homogeneous SIR models, in our model we use an infection rate that depends on the relative positions of the two parties. The infection rate between two positions will be given by the function $K$ defined on $D \times D$ with values in $ \R_+$. A susceptible individual $i$ becomes infected (in other words, his/her state changes from 0 to 1) at time t at rate (with some fixed $\gamma\in[0,1]$)
\begin{equation}
\label{1}
\begin{aligned}
\frac{1}{N^{1-\gamma}} \left[\displaystyle  \sum_{j\in\mathfrak{I}_N} \frac{K(X^{i},X^{j})}{\left[\displaystyle \sum_{\ell=1}^NK(X^{\ell},X^{j})\right]^{\gamma}} \lambda_{-j}(t)+ \sum_{j\in\mathfrak{S}_N} \frac{K(X^{i},X^{j})}{\left[\displaystyle \sum_{\ell =1}^NK(X^{\ell},X^{j})\right]^{\gamma}} \lambda_{j}(t- \tau_j^N)\right]\,,
\end{aligned}
\end{equation}
where $\tau_j^N$ is the infection time of the initially susceptible individual $j$ (in case individual $j$ never gets infected, then $\tau^N_j=+\infty$).
Denoting by $E^{i}_t$ the state of individual $i$ at time $t$, we have

\begin{equation}
\begin{aligned}
\mu_t^{S,N} &=\sum_{i\in\mathfrak{S}_N} \mathbf{1}_{E_{0}^{i}=0}\delta_{X^{i}}- \sum_{i\in\mathfrak{S}_N} \mathbf{1}_{t\geq \tau^N_i}\delta_{X^{i}}=\mu^{S,N}_0- \sum_{i\in\mathfrak{S}_N} \mathbf{1}_{t\geq \tau^N_i}\delta_{X^{i}},\\
\mu_t^{I,N}&=\sum_{i\in\mathfrak{I}_N} \mathbf{1}_{E_{t}^{i}=1}\delta_{X^{i}}+  \sum_{i\in\mathfrak{S}_N} \mathbf{1}_{E_{t}^{i}=1}\delta_{X^{i}}\\
&=\mu^{I,N}_0- \sum_{i\in\mathfrak{I}_N} \mathbf{1}_{\eta_{-i}\leq  t}\delta_{X^{i}}+\sum_{i\in\mathfrak{S}_N} \mathbf{1}_{t\leq \tau^N_i}\delta_{X^{i}}-
\sum_{i\in\mathfrak{S}_N} \mathbf{1}_{ t \geq \tau^N_i+ \eta_i  }\delta_{X^{i}},\\
\mu_t^{R,N}&=\sum_{i\in\mathfrak{R}_N} \mathbf{1}_{E_{0}^{i}=2}\delta_{X^{i}}+\sum_{i\in\mathfrak{I}_N} \mathbf{1}_{E_{t}^{i}=2}\delta_{X^{i}}+  \sum_{i\in\mathfrak{S}_N} \mathbf{1}_{E_{t}^{i}=2}\delta_{X^{i}}\\
&=\mu_0^{R,N}+\sum_{i\in\mathfrak{I}_N} \mathbf{1}_{\eta_{-i}\leq t}\delta_{X^{i}}+\sum_{i\in\mathfrak{S}_N} \mathbf{1}_{ \tau^N_i+ \eta_i \leq t }\delta_{X^{i}},\\
\mu_t^{N}&= \mu_t^{S,N}+\mu_t^{I,N}+\mu_t^{R,N}
\\
&=\sum_{i=1}^N\delta_{X^{i}}:=\mu^{N},\\
\mu_t^{\mathfrak{F},N} &=  \sum_{i\in\mathfrak{I}_N} \lambda_{-i}(t)\delta_{X^{i}} +  \sum_{i\in\mathfrak{S}_N} \lambda_{i}(t- \tau_i^N)\delta_{X^{i}},
\end{aligned}
\end{equation}
where 
\begin{itemize}
	\item $\mu_t^{S,N}$ is the empirical measure of susceptible individuals at time t;
	\item $\mu_t^{\mathfrak{F},N}$ is the empirical measure of the total force of infection at time t ;
	\item $\mu_t^{I,N}$ is the empirical measure of infected individuals at time t ;
	\item $\mu_t^{R,N}$ is the empirical measure of recovered individuals at time t;
	\item $\mu^{N}$ is the empirical measure of the total population, which does not depend upon $t$.
\end{itemize}
Now, we define $\overline{\mu}_t^{S,N}:=\dfrac{1}{N}\mu_t^{S,N}$; $\overline{\mu}_t^{I,N}:=\dfrac{1}{N}\mu_t^{I,N}$; $\overline{\mu}_t^{R,N}:=\dfrac{1}{N}\mu_t^{R,N}$; $\overline{\mu}^{N}:=\dfrac{1}{N}\mu^{N}$ and $\overline{\mu}_t^{\mathfrak{F},N}:=\dfrac{1}{N}\mu_t^{\mathfrak{F},N}$.
We rewrite \eqref{1} as follows. We first define
\[ M^N(x,y):=\frac{K(x,y)}{\left[\int_DK(z,y)\overline{\mu}^{N}(dz) \right]^{\gamma}},\]
so that the rate of infection of the individual $i$ at time $t$ can be written as
\begin{align*}
\mathbf{1}_{E_{t}^{i,N}=0} (\overline{\mu}_t^{\mathfrak{F},N},M^N(X^{i},\cdot))\,.
\end{align*}
We shall discuss the role of the parameter $\gamma\in[0,1]$ below, after the statement of our main result. 

Let $\eta_j$ be the random variable defined by $\eta_j := \sup\{t > 0,\;\lambda_{j}(t) > 0\} \quad \forall j \in  \mathbb{Z}^*.$
The two sequences of random variables $\{\eta_{-j} , j\geq 1\}$ and $\{\eta_{j} , j\geq 1\}$  are i.i.d and globally independent of each other. $F(t) := \mathbb{P}(\eta_1\leq t) \text{ and }F_0(t) := \mathbb{P}(\eta_{-1}\leq t)$ are the distribution functions of $\eta_j$ for $j\in \mathbb{Z}_+$ and for $j\in \mathbb{Z}_-$, respectively. For $  i \in \mathfrak{S}_N$, consider a counting process $A_i^N(t),$ which takes the value $0$ when individual $i$ is not yet infected at time t, and takes the value $1$ when the latter has been infected by time t. Thus, $\tau_i^N := \inf\{t > 0,\; A_i^N(t) = 1\}$. We define $A_i^N$ as follows :
\begin{align}\label{Ai}
A_i^N(t)=\int_{0}^{t}\int_{0}^{\infty}\mathbf{1}_{A^N_i(s^{-}) = 0}\mathbf{1}_{u\leq (\overline{\mu}_{s^-}^{\mathfrak{F},N},M^N(X^i,\cdot))}P^{i}(ds,du),
\end{align}
where the $\{P^{i},\;i \geq 1\}$ are standard Poisson random measures on $\mathbb{R}_{+}^2$ which are mutually independent, and globally independent of the 
$\{X^i,\ 1\le i\le N;\ \lambda_j,\ j\in\mathbb{Z}^\ast\}$.

The next proposition follows readily from our model.
\begin{propo}
	For all $\varphi \in C_b(D)$, $\left\{ \overline{\mu}_t^{S,N},\overline{\mu}_t^{\mathfrak{F},N},\overline{\mu}_t^{I,N}, \overline{\mu}_t^{R,N},\ t\ge0\right\}$ satisfies
	\begin{equation}
	\label{g}
	\left\lbrace
	\begin{aligned}
	(\overline{\mu}_t^{S,N},\varphi)&= (\overline{\mu}_0^{S,N},\varphi)
	-\frac{1}{N} \sum_{i\in\mathfrak{S}_N}\varphi(X^{i})A^N_i(t),\\
	(\overline{\mu}_t^{\mathfrak{F},N},\varphi) &= \frac{1}{N} \sum_{i\in\mathfrak{I}_N} \lambda_{-i}(t)\varphi({X^{i}}) +\frac{1}{N} \sum_{i\in\mathfrak{S}_N} \lambda_{i}(t- \tau_i^N)\varphi({X^{i}}),\\
	(\overline{\mu}_t^{I,N},\varphi)&= (\overline{\mu}_0^{I,N},\varphi)+\frac{1}{N} \sum_{i\in\mathfrak{S}_N}\varphi(X^{i})A^N_i(t)- \frac{1}{N}\sum_{i\in\mathfrak{I}_N}\varphi(X^{-i})\mathbf{1}_{\eta_{-i}\leq t}\\
	&- \frac{1}{N} \sum_{i\in\mathfrak{S}_N}\varphi(X^{i})\int_0^t\mathbf{1}_{\eta_i\leq t-s}dA^N_i(s),\\
	(\overline{\mu}_t^{R,N},\varphi)&= (\overline{\mu}_0^{R,N},\varphi)+ \frac{1}{N}\sum_{i\in\mathfrak{I}_N}\varphi(X^{-i})\mathbf{1}_{\eta_{-i}\leq t}+\frac{1}{N} \sum_{i\in\mathfrak{S}_N}\varphi(X^{i})\int_0^t\mathbf{1}_{\eta_i\leq t-s}dA^N_i(s)\,.
	\end{aligned}
	\right.
	\end{equation}	
\end{propo}
\section{Law of large numbers of measures }
\label{kk}
In this section,  we determine the limits of the empirical measures defined in section \ref{k} when the population size tends to infinity. Necessary intermediate results are established; they are summarized in lemmas and propositions. In what follows, we are given a probability measure $\overline{\mu}$ on $D$, with the density $\overline{\mu}(z)$.

\subsection{Assumptions and statement of the main result}
Let us first formulate our assumptions. The first one concerns the domain $D\subset\R^d$. We recall that $D$ is supposed to be compact. We need to make an assumption concerning the regularity of its boundary, namely we assume the following interior cone condition. In order to formulate that assumption, we first define  the following cone: for any $y\in\R^d$, $\ell\in\R^d$ with $\|\ell\|=1$ and $0<\alpha <1$,
\[ C(y,\ell,\alpha)=\{z\in\R^d, 0\le(z-y,\ell)\le\alpha\|z-y\|\}\,.\]
Note that $\|\cdot\|$ stands for the norm defined by $\|z\|=\sqrt{\sum_{i=1}^dz_i^2}$. Our assumptions are the following. 
\begin{assumption}
	\label{A1}
	There exists $\alpha\in(0,1)$ and $r>0$ such that for any $y\in D$, there exists a unit vector $\ell_y$ such that
	\[ C(y,\ell_y,\alpha)\cap B(y,r)\subset D,\]
	where $B(y,r)$ denotes the ball of radius $r$ centered at $y$.
	\end{assumption}
	Note that the above assumption is in fact essentially an assumption about the boundary of $D$. It says that
	from any point $y\in\partial D$, we can choose a direction $\ell_y$ such that all points of $C(y,\ell_y,\alpha)\cap B(y,r)$
	are in $D$. As a result, the Lebesgue measure of $B(y,r)\cap D$ is bounded from below by a fixed constant, for any 
	$y\in D$, which will be crucial for us.
	
	The second assumption concerns the kernel $K$, the measure $\overline{\mu}$, and the initial condition.
	The same $r$ appears in both assumptions \ref{A1} and \ref{A2}.
\begin{assumption}
	\label{A2}
	We assume that:
	\begin{itemize}
		\item The function K is a measurable map from $D^2$ into $\R_+$,  s. t. $K(x,x)=0$; and  there exist $\underline{c} > 0$ s. t. for any $0<\|x-y\|\leq r$,  then $K(x,y)\geq \underline{c}$, while $K$ is bounded on $\{(x,y)\in D\times D,\ \|x-y\|>r\}$ ;
		\item $\displaystyle \sup_{x\in D} \int_D K(x,y)dy <\infty$, $\displaystyle \sup_{y\in D} \int_D K(x,y)dx <\infty$;
          \item $\displaystyle \int_D \int_D K^2(x,y)dxdy <\infty$;
		\item $ \displaystyle 0<\inf_{x\in D}\overline{\mu}(x)<\sup_{x\in D} \overline{\mu}(x)< \infty$; 
		\item $\overline{S}^N(0):=\frac{S^N(0)}{N}\rightarrow\overline{S}(0)$; $\overline{I}^N(0):=\frac{I^N(0)}{N}\rightarrow\overline{I}(0)$; and $\overline{R}^N(0):=\frac{R^N(0)}{N}\rightarrow\overline{R}(0)$ a.s. as $N\to\infty$;
		\item $\displaystyle \P(E^i_0=0)=\overline{S}(0),\quad \P(E^i_0=1)=\overline{I}(0), \quad \text{and}\quad \P(E^i_0=2)=\overline{R}(0)$;
\item The pairs $\left\{\left(E^{i}_0, X^i\right), i=1;\cdots;N\right\}$ are i.i.d;
\item For all $i=1;\cdots;N$, $X^i=\begin{cases}
X_S^i\quad \text{if}\quad E^i_0=0,\\
X_I^i\quad \text{if}\quad E^i_0=1,\\
X_R^i\quad \text{if}\quad E^i_0=2;
\end{cases}$
\item $\left(X_S^i, i=1;\cdots;N\right)$ are  i.i.d with the density function $\pi_S$, $\big(X_I^i, i=1;$ $\cdots;N\big)$ are  i.i.d with the density function $\pi_I$ and $\left(X_R^i, i=1;\cdots;N\right)$  are i.i.d with the density function $\pi_R$. Those three collections of r.v. are mutually independent.
			\end{itemize}
\end{assumption}
\begin{example} Let $\psi$ be a bounded measurable function from $D\times D$ into $\R_+$ and 
$\beta<d/2$. The following function $K$ defined by:
\[
K(x,y) =
\begin{cases}
0, & \text{if } x = y, \\[4pt]
\dfrac{1}{\|x - y\|^{\beta}}, & \text{if } 0<\|x - y\| \le r,  \\[6pt]
\psi(x,y), & \text{otherwise,}
\end{cases}
\]
satisfies the above assumptions.
\end{example}

\medskip

Of course, $\overline{\mu}(x)=\overline{S}(0)\pi_S(x)+\overline{I}(0)\pi_I(x)+\overline{R}(0)\pi_R(x)$ and $\overline{\mu}$ is a density function. For an epidemic to develop, we need that both $\overline{S}(0),\overline{I}(0)>0$.

We shall use the notation $X^i$  to denote the position of the individual $i$, and the notation 
$X_S$ (resp. $X_I$, $X_R$) to denote a position following the distribution
$\pi_S$ (resp. $\pi_I$, $\pi_R$).

We next define the notation
\[ M(x,y):=\frac{K(x,y)}{\left[\int_DK(z,y)\overline{\mu}(dz) \right]^{\gamma}}\,.\]

We can now state our main result.
\begin{theorem}
	\label{theo1}
	Under assumptions \ref{A1} and \ref{A2}, the sequence $(\overline{\mu}^{S,N}, \overline{\mu}^{\mathfrak{F},N},\overline{\mu}^{I,N}$,
	$\overline{\mu}^{R,N})_{N \geq 1}$ converges in probability in  $\mathbb{D}^4_{\mathcal{M}}$ to $(\overline{\mu}^{S}, \overline{\mu}^{\mathfrak{F}},\overline{\mu}^{I},\overline{\mu}^{R})$ such that for all $\varphi \in C_b(D)$, $\left\{ (\overline{\mu}^S_t,\varphi), (\overline{\mu}^{\mathfrak{F}}_t,\varphi),(\overline{\mu}^I_t,\varphi),(\overline{\mu}^R_t,\varphi),   t\ge0 \right\}$  satisfies
	\begin{equation}
	\label{sy}
	\left\lbrace
	\begin{aligned}
	(\overline{\mu}_t^{S},\varphi)&= (\overline{\mu}_0^{S},\varphi)- \int_0^t\int_{D}\varphi(x)(\overline{\mu}_s^{\mathfrak{F}},M(x,\cdot)\overline{\mu}_s^{S}(dx)ds,\\
	(\overline{\mu}_t^{\mathfrak{F}},\varphi) &=\overline{\lambda}^0(t)(\overline{\mu}^I_0,\varphi)+\int_0^t\overline{\lambda}(t-s)\int_{D}\varphi(x)(\overline{\mu}_s^{\mathfrak{F}},M(x,\cdot)\overline{\mu}_s^{S}(dx)ds,\\
	(\overline{\mu}_t^{I},\varphi)&= (\overline{\mu}_0^{I},\varphi)F_0^c(t)+\int_0^tF^c(t-s)\int_{D} \varphi (x)(\overline{\mu}_s^{\mathfrak{F}},M(x,\cdot)\overline{\mu}_s^S(dx)ds,\\
	(\overline{\mu}_t^{R},\varphi)&=(\overline{\mu}_0^{R},\varphi)+ (\overline{\mu}_0^{I},\varphi)F_0(t)+\int_0^t F(t-s) \int_{D} \varphi (x)(\overline{\mu}_s^{\mathfrak{F}},M(x,\cdot)\overline{\mu}_s^S(dx)ds,\\
	\overline{\mu}_0^{S}(dx)&= \overline{S}(0)\pi_S(x)dx;\ \overline{\mu}_0^{I}(dx)= \overline{I}(0)\pi_I(x)dx;\ \overline{\mu}_0^{R}(dx)= \overline{R}(0)\pi_R(x)dx, \\
	 \overline{\mu}(dx)&=\overline{\mu}^S_t(dx)+\overline{\mu}^I_t(dx)+\overline{\mu}^R_t(dx),\forall t\ge0.
	\end{aligned}
	\right.
	\end{equation}	
\end{theorem}

We now want to discuss the role of our parameter $\gamma\in[0,1]$.
Suppose first that the function $y\mapsto \int_DK(z,y)\overline{\mu}(dz)$ is constant. We can renormalize the kernel $K$ in such a way that this constant is $1$. Then $M(x,y)$ does not vary with $\gamma$, hence also for large $N$, $M^N(x,y)$ does not vary much with $\gamma$. We then assume that 
$y\mapsto \int_DK(z,y)\overline{\mu}(dz)$ is not constant. Assume for convenience that $K$ has been renormalized in such a way that $\displaystyle \int_D\int_DK(z,y)\overline{\mu}(dz)\overline{\mu}(dy)=1$. Hence for large $N$, $\displaystyle \int_DK(z,y)\overline{\mu}^N(dz)$
fluctuates around $1$. In an homogeneous model, $\overline{S}^N_t=S^N_t/N$ can be thought of as the probability that an individual met at random is susceptible. In our spatial model,
in the case $\gamma=1$, we may think that $\displaystyle \int_DM^N(x,y)\overline{\mu}^{S,N}_t(dx)$ is the probability that the individual met by an infected individual located at $y$ at time $t$ is susceptible.
The values $\gamma<1$ tend to increase the rate of infection of individuals located in more populated regions (which is likely to correspond to reality), the effect being 
more significant with smaller $\gamma$. 

\subsection{Study of the system of limiting equations}
We will first formulate a basic result concerning the system \eqref{sy}.
\begin{propo}
	\label{oo}
	The system of equations \eqref{sy} admits a unique solution  $\left\{\left(\overline{\mu}^S_t,\overline{\mu}_t^{\mathfrak{F}},\overline{\mu}_t^{I},\overline{\mu}_t^{R}\right),\ t\ge0\right\}$   which is absolutely continuous with respect to the Lebesgue measure, with the densities\\ $\left\{\left(\overline{\mu}^S(t,.), \overline{\mu}^{\mathfrak{F}}(t,.), \overline{\mu}^{I}(t,.) ,\overline{\mu}^{R}(t,.)\right),\ t\ge0\right\}$ satisfying for all $(t,x)\in \R_+\times D$
	\begin{equation}
	\label{syst}
	\left\lbrace
	\begin{aligned}
	\overline{\mu}^{S}(t,x)&= \overline{\mu}^{S}(0,x)- \int_0^t\overline{\mu}^S(s,x)\int_DM(x,y)\overline{\mu}^{\mathfrak{F}}(s,y)dyds,\\
	\overline{\mu}^{\mathfrak{F}}(t,x) &=\overline{\lambda}^0(t)\mu^I(0,x)+\int_0^t\overline{\lambda}(t-s)\overline{\mu}^S(s,x)\int_DM(x,y)\overline{\mu}^{\mathfrak{F}}(s,y)dyds,\\
	\overline{\mu}^{I}(t,x)&= \mu^I(0,x)F_0^c(t)+\int_0^t F^c(t-s)\overline{\mu}^S(s,x)\int_DM(x,y)\overline{\mu}^{\mathfrak{F}}(s,y)dyds,\\
	\overline{\mu}^{R}(t,x)&= \overline{\mu}^R(0,x)+\overline{\mu}^I(0,x)F_0(t)+\int_0^tF(t-s) \overline{\mu}^S(s,x)\int_DM(x,y)\overline{\mu}^{\mathfrak{F}}(s,y)dyds,\\
	\overline{\mu}^{S}(0,x)&= \overline{S}(0)\pi_S(x),\ \overline{\mu}^{I}(0,x)= \overline{I}(0)\pi_I(x),\ \overline{\mu}^{R}(0,x)= \overline{R}(0)\pi_R(x),
	\\ \overline{\mu}(x)&=\overline{\mu}^S(t,x)+\overline{\mu}^{I}(t,x)+\overline{\mu}^{R}(t,x),\ \forall t\ge0.
	\end{aligned}
	\right.
	\end{equation}	
\end{propo}
Admitting for a moment the first part of Proposition \ref{oo}, we first establish the following a priori estimates.
\begin{propo}\label{pro:aprioriestim}
	Let $T>0$, and let $(\overline{\mu}^S,\overline{\mu}^{\mathfrak{F}} )$  be a solution of the first two equations of \eqref{syst}. Then there exists  positive constants $C$ 	
	and $c$ such that:
	\begin{itemize}
		\item $\forall t\in[0;T],$ $\|\overline{\mu}^{S}(t,.)\|_{\infty}\leq C;$
		\item $\displaystyle \inf_{y\in D} \int_{D}K(z,y)\overline{\mu}(z)dz \ge c$;
		\item $\forall t\in[0;T]$; $\|\overline{\mu}^{\mathfrak{F}}(t,.)\|_{\infty}\leq C.$
	\end{itemize}
\end{propo}
\begin{proof}
For any $(t,x)\in [0;T]\times D$, 
\begin{align}
\overline{\mu}^{S}(t,x)&\leq \overline{\mu}^{S}(0,x)\leq \overline{\mu}(x)\nonumber\\
\|\overline{\mu}^{S}(t,.)\|_{\infty}&\leq\|\overline{\mu}(.)\|_{\infty}:=C.\\
\int_{D}K(z,y)\overline{\mu}(dz)&= \int_{D}K(z,y)\overline{\mu}(z)dz\geq \inf_{z\in D} \overline{\mu}(z) \int_{D}K(z,y)dz\nonumber\\
&\geq \inf_{z\in D} \overline{\mu}(z) \int_{D\cap B(y,r)}K(z,y)dz\nonumber\\
&\geq \inf_{z\in D} \overline{\mu}(z)\,\underline{c}\int_{C(y,\ell_y,\alpha)\cap B(y,r)}dz=c.
\end{align}
Next
\begin{align}
\overline{\mu}^{\mathfrak{F}}(t,x) &=\overline{\lambda}^0(t)\mu^I(0,x)+\int_0^t\overline{\lambda}(t-s)\overline{\mu}^S(s,x)\int_DM(x,y)\overline{\mu}^{\mathfrak{F}}(s,y)dyds,\nonumber\\
\overline{\mu}^{\mathfrak{F}}(t,x) &\leq \lambda^*\overline \mu(x)+\lambda^*\int_0^t\overline{\mu}^S(s,x)\int_DM(x,y)\overline{\mu}^{\mathfrak{F}}(s,y)dyds,\nonumber\\
\label{p1}
\|\overline{\mu}^{\mathfrak{F}}(t,.)\|_{\infty}&\leq \lambda^*C+\lambda^*\int_0^t\|\overline{\mu}^{S}(s,.)\int_DM(\cdot,y)\overline{\mu}^{\mathfrak{F}}(s,y)dy\|_{\infty}ds,\\
\overline{\mu}^{S}(t,x)\int_D&M(x,y)\overline{\mu}^{\mathfrak{F}}(t,y)dy=\int_{D} \dfrac{K(x,y)\overline{\mu}^{S}(t,x)}{\left[\displaystyle \int_{D}K(z,y)\overline{\mu}(z)dz \right]^{\gamma}}\overline{\mu}^{\mathfrak{F}}(t,y)dy\nonumber\\
&\quad\quad\quad\quad\quad\quad\leq \frac{C}{c^\gamma} \|\overline{\mu}^{\mathfrak{F}}(t,.)\|_{\infty}\int_{D} K(x,y)dy\,.\nonumber
\end{align}
By assumption, $\displaystyle \sup_{x\in D}\int_{D} K(x,y)dy < \infty$. Thus,
\begin{align}
\label{p2}
\|\overline{\mu}^{S}(t,.)\int_DM(\cdot,y)\overline{\mu}^{\mathfrak{F}}(t,y)dy\|_{\infty}&\leq C\|\overline{\mu}^{\mathfrak{F}}(t,.)\|_{\infty}\,.
\end{align}
	From \eqref{p1} and \eqref{p2}, we deduce that  $\forall t\in[0;T]$ 
\begin{align*}
\|\overline{\mu}^{\mathfrak{F}}(t,.)\|_{\infty}&\leq \lambda^*C+\lambda^*C\int_0^t\|\overline{\mu}^{\mathfrak{F}}(s,.)\|_{\infty}ds,
\end{align*}
which combined with Gronwall's inequality yields
\begin{align}
\|\overline{\mu}^{\mathfrak{F}}(t,.)\|_{\infty}&\leq\lambda^*Ce^{\lambda^*CT},\quad  \forall t\in [0;T]\,.
\end{align}
\end{proof}
It follows in particular from Proposition \ref{pro:aprioriestim} that $M(x,y)\le c^{-\gamma} K(x,y)$.
\begin{proof}\textit{of Proposition \ref{oo}}.
We first show that for all $t\ge0$ any solution $\left(\overline{\mu}^S_t,\overline{\mu}_t^{\mathfrak{F}},\overline{\mu}_t^{I},\overline{\mu}_t^{R}\right)$ of \eqref{sy}  is absolutely continuous with respect to the Lebesgue measure, and the densities $\big(\overline{\mu}^S(t,.), \overline{\mu}^{\mathfrak{F}}(t,.), \overline{\mu}^{I}(t,.),$ $\overline{\mu}^{R}(t,.)\big)$ verify \eqref{syst}.\\
 From the first equation of \eqref{sy}, $\overline{\mu}^S_t\leq \overline{\mu}^S_0$. Since $\overline{\mu}^S_0$ is absolutely continuous,  $\overline{\mu}^S_t$ has the same property, and we denote its density by $\overline{\mu}^S(t,x)$. \\
From the third equation of \eqref{sy}, $\displaystyle \overline{\mu}_t^{I}\leq  \overline{\mu}_0^{I}+ \int_0^t \overline{\Gamma}(s,.)\overline{\mu}^{S}_sds$, thus $\overline{\mu}_t^{I}$ is absolutely continuous , since $\overline{\mu}^I_0$ is absolutely continuous, as well as $\overline{\mu}^{S}_s$ for all s. 
The same argument applies  to $ \overline{\mu}_t^{\mathfrak{F}}$ and $ \overline{\mu}_t^{R}$.The system of equation \eqref{syst} now follows readily from \eqref{sy}.\\
We will verify that  $\left(\overline{\mu}^S(t,.), \overline{\mu}^{\mathfrak{F}}(t,.), \overline{\mu}^{I}(t,.) ,\overline{\mu}^{R}(t,.)\right)$ is unique.  For that sake, it suffices to show that the solution $(\overline{\mu}^{S}(t,.), \overline{\mu}^{\mathfrak{F}}(t,.))$ of the first two equations of the system is unique.
	The first two equations of the system \eqref{syst} constitute the following system
	\begin{equation}
	\label{sys}
	\left\lbrace
	\begin{aligned}
	\overline{\mu}^{S}(t,x)&= \overline{\mu}^{S}(0,x)- \int_0^t\int_{D}M(x,y)\overline{\mu}^{\mathfrak{F}}(s,y)\overline{\mu}^S(s,x)dyds,\nonumber\\
	\overline{\mu}^{\mathfrak{F}}(t,x) &=\overline{\lambda}^0(t)\mu^I(0,x)+\int_0^t\overline{\lambda}(t-s)\int_{D}M(x,y)\overline{\mu}^{\mathfrak{F}}(s,y)\overline{\mu}^{S}(s,x)dyds,\nonumber\\
 \overline{\mu}^{S}(0,x)&= \overline{S}(0)\pi_S(x); \overline{\mu}^{I}(0,x)= \overline{I}(0)\pi_I(x), \text{and}\quad   \overline{\mu}^{R}(0,x)= \overline{R}(0)\pi_{R}(x),\\ \overline{\mu}(x)&=\overline{\mu}^S(0,x)+ \overline{\mu}^I(0,x)+\overline{\mu}^R(0,x)\,.\nonumber
	\end{aligned}
	\right.
	\end{equation}
	Note that this is a sytem of ODE - integral equation, which depends upon  the parameter $x\in D$. The pair of equations for the various values of $x$ are coupled through the integrals over $D$.
	
	Let $\left(f_1(t,.), g_1(t,.)\right)$ and $\left(f_2(t,.), g_2(t,.)\right)$ be two solutions of the above system with the same initial condition.\\
	On the one hand, exploiting Proposition \ref{pro:aprioriestim}, we obtain	
	\begin{align}
	f_1(t,x)-f_2(t,x)&= \int_0^t  \left(f_2(s,x)-f_1(s,x)\right)\int_{D}M(x,y)g_2(s,y)dyds\nonumber\\
	&+\int_0^tf_1(s,x)\int_{D}M(x,y)\left( g_2(s,y)-g_1(s,y)\right)dyds,\nonumber\\
	\|f_1(t,.)-f_2(t,.)\|_{\infty}&\leq \frac{C}{c^\gamma}\int_0^t \|f_2(s,.)-f_1(s,.)\|_{\infty}\|g_2(s,.)\|_\infty ds\nonumber\\
	&+\frac{C}{c^\gamma}\int_0^t\| g_2(s,.)-g_1(s,.)\|_{\infty} \|f_1(s,.)\|_{\infty}ds,\nonumber\\
	\label{C}
	\|f_1(t,.)-f_2(t,.)\|_{\infty}&\leq C\int_0^t \left(\|f_2(s,.)-f_1(s,.)\|_{\infty}+\| g_2(s,.)-g_1(s,.)\|_{\infty}\right) ds\,.
	\end{align}
	Moreover,
	\begin{align}
	g_1(t,x)-g_2(t,x)&=\int_0^t\overline{\lambda}(t-s) (f_1(s,x)-f_2(s,x))\int_{D}M(x,y)g_1(s,y)dyds\nonumber\\
	&+\int_0^t \overline{\lambda}(t-s)f_2(s,x)\int_{D} (g_1(s,y)-g_2(s,y))M(x,y)dyds,\nonumber\\
	|g_1(t,x)-g_2(t,x)|&\leq C\int_0^t \| f_1(s,.)-f_2(s,.)\|_{\infty}\left\|\int_{D}M(\cdot,y)dy\right\|_{\infty}ds\nonumber\\
	&+C \int_0^t \|g_1(s,.)-g_2(s,.)\|_{\infty}\left\|\int_{D} M(\cdot,y)dy\right\|_{\infty}ds,\nonumber\\
	\label{C1}
\|g_1(t,.)-g_2(t,.)\|_{\infty}	&\leq C\int_0^t \left(\| f_1(s,.)-f_2(s,.)\|_{\infty}+\|g_1(s,.)-g_2(s,.)\|_{\infty}\right)ds\,.
	\end{align}
	From \eqref{C} and \eqref{C1}, we have
	\begin{align*}
	\|g_1(t,.)-g_2(t,.)\|_{\infty}&+\|f_1(t,.)-f_2(t,.)\|_{\infty}\\&\leq C \int_0^t  \left(\|g_1(s,.)-g_2(s,.)\|_{\infty}+\|f_1(s,.)-f_2(s,.)\|_{\infty}\right)ds\,.
	\end{align*}
	Using Gronwall's inequality, we obtain 
	\begin{align*}
	\|g_1(t,.)-g_2(t,.)\|_{\infty}+\|f_1(t,.)-f_2(t,.)\|_{\infty}&=0.
	\end{align*}
	We finally stress that since any solution to \eqref{sy} has a density which solves \eqref{syst}, uniqueness of \eqref{syst} implies uniqueness of \eqref{sy}. On the other hand, we could prove existence of a solution to \eqref{syst}, but this will follow from Theorem \ref{theo1}, which provides existence of a solution to \eqref{sy}, hence of a solution to \eqref{syst}.
\end{proof}

\subsection{A variant of the sequence indexed by $N$}
The lower bound in Proposition \ref{pro:aprioriestim} tells us that $\inf_{y\in D}\int_D K(z,y)\overline{\mu}(dz)\ge c$, hence the denominator in the expression of $M(x,y)$ is lower bounded by $c^\gamma$. However, in the $N$ model
$\overline{\mu}$ is replaced by $\overline{\mu}^N$, and $\inf_{y\in D}\int_D K(z,y)\overline{\mu}^N(dz)$ is not lower bounded
a.s. for all $N$, which clearly creates a difficulty for our proof. However, it is a sequence whose limit in probability is $\ge c$, hence the probability that it is $\ge c/2$ 
converges to $1$ as $N\to\infty$. This is the motivation for modifying our model, by replacing the denominator by a quantity which is uniformly bounded away from $0$ as follows.
Let $\Phi$ be the continuous function from $\mathbb{R}_+$ into $\mathbb{R}_+$ defined by $\Phi(u):=\left( u\vee \delta\right)^\gamma$, where $\delta<c/2$ will be specified in the last proof below. We easily verify that 
$\left|\left(\frac{1}{\Phi}\right)'(u)\right|\le\gamma\delta^{-\gamma-1}=:k_\Phi$ As a consequence, for any $u,v>0$,
\begin{equation}\label{LipPhi}
\left|\frac{1}{\Phi(u)}-\frac{1}{\Phi(v)}\right|\le k_\Phi |u-v|\,.
\end{equation}

We next define
\begin{align}
\M^N(x,y)&:=\frac{K(x,y)}{\Phi\left(\int_DK(z,y)\overline{\mu}^N(dz)\right)},\nonumber\\
\tilde{A}^N_i(t)&:= \int_0^t \int_0^{\infty}\mathbf{1}_{\tilde{A}^N_i(s^{-})=0}\mathbf{1}_{u \leq (\tilde{\mu}_{s^-}^{\mathfrak{F},N},\M^N(X^i,\cdot)) }P^{i}(ds,du),\label{Aitilde}
\end{align}
where $\tilde{\mu}_{t}^{\mathfrak{F},N}$ is defined below.

The variant $\left(\tilde{\mu}^{S,N}, \tilde{\mu}^{I,N},\tilde{\mu}^{R,N},\tilde{\mu}^{\mathfrak{F},N} \right)$ of $\left(\overline{\mu}^{S,N}, \overline{\mu}^{I,N},\overline{\mu}^{R,N},\overline{\mu}^{\mathfrak{F},N} \right)$ verifies, for all $t\in [0;T]$, $\varphi \in C_b\left(D\right)$
\begin{align*}
(\tilde{\mu}_t^{S,N},\varphi)&= (\overline{\mu}_0^{S,N},\varphi)-\frac{1}{N} \sum_{i\in\mathfrak{S}_N}\varphi(X^{i})\tilde{A}^N_i(t)\\
(\tilde{\mu}_t^{I,N},\varphi)&= (\overline{\mu}_0^{I,N},\varphi)+\frac{1}{N} \sum_{i\in\mathfrak{S}_N}\varphi(X^{i})\tilde{A}^N_i(t)- \frac{1}{N}\sum_{i\in\mathfrak{I}_N}\varphi(X^{i})\mathbf{1}_{\eta_{-i}\leq t}\\
&\quad- \frac{1}{N} \sum_{i\in\mathfrak{S}_N}\varphi(X^{i})\int_0^t\mathbf{1}_{\eta_i\leq t-s}d\tilde{A}^N_i(s)\\
(\tilde{\mu}_t^{R,N},\varphi)&= (\overline{\mu}_0^{R,N},\varphi)+ \frac{1}{N}\sum_{i\in\mathfrak{I}_N}\varphi(X^{i})\mathbf{1}_{\eta_{-i}\leq t}+\frac{1}{N} \sum_{i\in\mathfrak{S}_N}\varphi(X^{i})\int_0^t\mathbf{1}_{\eta_i\leq t-s}d\tilde{A}^N_i(s)\\
(\tilde{\mu}_t^{\mathfrak{F},N},\varphi)&=  \frac{1}{N}\sum_{i\in\mathfrak{I}_N} \lambda_{-i}(t)\varphi(X^{i})+\frac{1}{N}  \sum_{i\in\mathfrak{S}_N} \lambda_{i}(t- \tilde{\tau}_i^N)\varphi(X^{i}),\\
 \tilde{\tau}_i^N&:=\inf \{ t\geq 0;\tilde{A}^N_i(t)=1\}
\end{align*}

Before proving our main result Theorem \ref{theo1}, we shall prove the following result.
\begin{theorem}
	\label{th}
	Under assumptions \ref{A1} and \ref{A2}, 
	the sequence $(\tilde{\mu}^{S,N}, \tilde{\mu}^{\mathfrak{F},N},\tilde{\mu}^{I,N}$, $\tilde{\mu}^{R,N})_{N \geq 1}$ converges in probability in  $\mathbb{D}_{\mathcal{M}}^4$ to $\left(\overline{\mu}^S,\overline{\mu}^{\mathfrak{F}},\overline{\mu}^{I}, \overline{\mu}^{R}\right)$, the unique solution of \eqref{sy}.
	\end{theorem}
We note that since the denominator in the expression for $M(x,y)$ is lower bounded by $c$, if we define 
$\M(x,y):=\frac{K(x,y)}{\Phi\left(\displaystyle \int_DK(z,y)\overline{\mu}(dz)\right)}$, we have that $\M(x,y)=M(x,y)$. 
As a consequence, $\left(\tilde{\mu}^S,\tilde{\mu}^{\mathfrak{F}},\tilde{\mu}^{I},\tilde{\mu}^{R}\right)=\left(\overline{\mu}^S,\overline{\mu}^{\mathfrak{F}},\overline{\mu}^{I},\overline{\mu}^{R}\right)$, where the quantities with $\tilde{}$ are defined as for the $N$ model by replacing $M$ by $\M$.

\subsection{A McKean-Vlasov equation}

Recall that in order to find the solution of \eqref{sy}, it suffices to find a solution of the system \eqref{syst}. Moreover, 
it suffices to solve the first two equations of \eqref{syst} for $\left(\overline{\mu}^S,\overline{\mu}^{\mathfrak{F}}\right)$.
In fact,  $\overline{\mu}^S$ solves the first equation iff (from now on, we write $(\overline{\mu}^{\mathfrak{F}}_t,M(x,\cdot))$
for $\displaystyle \int_DM(x,y)\overline{\mu}^{\mathfrak{F}}(t,y)dy$)
\begin{align}\label{explicit}
\overline{\mu}^S(t,x)=\overline{\mu}^S(0,x)\exp\left(-\int_0^t(\overline{\mu}^{\mathfrak{F}}_s,M(x,\cdot))ds\right)\,.
\end{align}
Consequently $\left(\overline{\mu}^S,\overline{\mu}^{\mathfrak{F}}\right)$ solve the first two equations of \eqref{syst}
iff $\overline{\mu}^{\mathfrak{F}}$ solves
\begin{align}\label{eq:F}
&\overline{\mu}^{\mathfrak{F}}(t,x)=\overline{\lambda}^0(t)\mu^I(0,x)\nonumber\\&+\int_0^t\overline{\lambda}(t-s)
\overline{\mu}^S(0,x)\exp\left(-\int_0^s(\overline{\mu}^{\mathfrak{F}}_r,M(x,\cdot))dr\right)(\overline{\mu}^{\mathfrak{F}}_s,M(x,\cdot))ds\,,
\end{align}
and $\overline{\mu}^S$ is given by \eqref{explicit}. An alternative proof of Proposition \ref{oo} would have consisted to prove existence and uniqueness of the solution of \eqref{eq:F}, which would have been as easy.

Now, given $P$  a standard Poisson random measure on $\mathbb{R}_{+}^2$, we define the process 
$\{A(t,x),\ t\ge0,\ x\in D\}$ as follows
\begin{align}
\label{le}
A(t,x)&=\int_0^t\int_0^{\infty} \mathbf{1}_{A(s^{-},x)=0}\mathbf{1}_{u \leq ({F}(s),M(x,\cdot))}P(ds,du),\ t\ge0\ \text{ with}\nonumber\\
{F}(t,y)   &:= \overline{I}(0)\overline{\lambda}^0(t)\pi_I(y)
+\overline{S}(0)\E[\overline{\lambda}(t-{\tau}(y))]\pi_S(y)\\
\text{where}\quad{\tau}(x)&:=  \inf \{ t> 0 ;   A(t,x)=1  \},\  ({F}(t),M(x,\cdot))=\int_DM(x,y){F}(t,y)dy\nonumber.
\end{align}
Note that equation \eqref{le} is an equation of the McKean--Vlasov type, in the sense that the coefficient of that equation
depends upon the law of its solution. 
Indeed ${F}$ depends upon the laws of $\{{\tau}(y), y\in D\}$, that is upon
the laws of the processes $\{A(\cdot,y),\ y\in D\}$. 
Note that McKean-Vlasov type of equations are related to 
the ``propagation of chaos'', see \cite{szni}. We now establish existence and uniqueness of the solution to \eqref{le}.
The result is an extension to spatial models of a result in \cite{b}. However the argument is new and simpler than in \cite{b}. 
\begin{lemma}
\label{lll}
Equation \eqref{le} has a unique solution $\{A(t,x),\ t\geq 0,\ x\in D\}$, which is such that  ${F}(t,x)=\overline{\mu}^{\mathfrak{F}}(t,x)$, the unique solution to \eqref{eq:F}.	
\end{lemma}
\begin{proof}
We first note that if $A(t,x)$ solves \eqref{le}, then, since $\overline{S}(0)\pi_S(x)=\overline{\mu}^S(0,x)$
and $\overline{\lambda}(t)=0$ for $t<0$,
\begin{align*}
\overline{S}(0)&\E[\overline{\lambda}(t-{\tau}(x))]\pi_S(x)\\
&=
\overline{\mu}^S(0,x)\E\int_0^t\int_0^\infty\overline{\lambda}(t-s)\mathbf{1}_{A(s^{-},x)=0}\mathbf{1}_{u \leq ({F}(s),M(x,\cdot))}P(ds,du)\\
&=\overline{\mu}^S(0,x)\int_0^t\overline{\lambda}(t-s)\P(A(s,x)=0)({F}(s),M(x,\cdot))ds\\
&=\overline{\mu}^S(0,x)\int_0^t\overline{\lambda}(t-s)\exp\left(-\int_0^s({F}(r),M(x,\cdot))dr\right)({F}(s),M(x,\cdot))ds\\
&=\int_0^t\overline{\lambda}(t-s)({F}(s),M(x,\cdot))\overline{\mu}^S(s,x)ds\,,
\end{align*}
where we have used \eqref{explicit} for the last identity.
Hence the second line of \eqref{le} tells us that necessarily ${F}(t,x)$ solves \eqref{eq:F}. Since the latter has a unique solution, ${F}(t,x)$ is unique, hence $A(t,x)$ is unique as well. But this implies also existence: the first line
of \eqref{le} with  ${F}(t,x)$ replaced by $\overline{\mu}^{\mathfrak{F}}(t,x)$, the unique solution of \eqref{eq:F}, defines
a process $A(t,x)$, which, thanks to the above computation, solves \eqref{le}.
\end{proof}

For each $ i \in \mathfrak{S}_N$, we denote by $A_i(t,x)$ the solution of \eqref{le} with $P$ replaced by $P^i$, the PRM which appears in the definition of both $A^N_i(t)$ and $\tilde{A}^N_i(t)$, see \eqref{Ai} and \eqref{Aitilde}. We next define $A_i(t):=A_i(t,X^i)$, hence 
$\tau_i=\tau_i(X^i)$, and we remark that the $\{\tilde{A}^N_i(\cdot)-A_i(\cdot),\ i\in \mathfrak{S}_N\}$ are identically distributed. 
Note that we shall use below $\tau(X_S)$, which depends upon the two mutually independent random inputs
 $P$ and $X_S$.
 
 Clearly the $A_i$'s are i.i.d., which is not the case of the $A^N_i$'s or of the $\tilde{A}^N_i$'s. The fact that the $A_i$'s are i.i.d. will allow us to use the law of large numbers in order to take the limit as $N\to\infty$ in $N^{-1}\sum_{i\in\mathfrak{S}_N}A_i(t)$, and other quantities related to this.  The next Lemma will be crucial, in order to conclude concerning $N^{-1}\sum_{i\in\mathfrak{S}_N}\tilde{A}^N_i$.
\begin{lemma}
	\label{Col}
	For all $T>0$,
	$\displaystyle\sup_{i\in\mathfrak{S}_N} \mathbb{E}\left(\sup_{0 \leq t \leq T}\left| \tilde{A}^N_i(t)-A_i(t)\right|\right) \xrightarrow [N \rightarrow \infty]{} 0.$
\end{lemma}
\begin{proof} In this proof we shall use the following notations
\[\G(t,x):=(\overline{\mu}^{\mathfrak{F}}_t,M(x,\cdot)),\quad \GN(t,x):=(\tilde{\mu}^{\mathfrak{F},N}_t,\M^N(x,\cdot))\,.\]
	For any $t\in[0;T]$, $i\in\mathfrak{S}_N$, we have
	\begin{align}
	\left|\tilde{A}^N_i(t)- A_i(t)\right|&\leq \int_0^t \int_{\GN(s,X^i) \wedge\G(s,X^i)}^{\GN(s,X^i) \vee\G(s,X^i)} P^i(ds,du),
	\nonumber\\
	\sup_{0 \leq t \leq T}\left|\tilde{A}^N_i(t)- A_i(t)\right|&\leq \int_0^T \int_{\tilde{\Gamma}^N(s,X^i) \wedge
	\G(s,X^i)}^{\tilde{\Gamma}^N(s,X^i) \vee\G(s,X^i)} P^i(ds,du),\nonumber\\
	\label{t2}
	\mathbb{E}\left(\sup_{0 \leq t \leq T}\left|\tilde{A}^N_i(t)- A_i(t)\right|\right)&\leq \int_0^T \mathbb{E}\left(\left|\tilde{\Gamma}^N(s,X^i)-\G(s,X^i)\right|\right)ds.
	\end{align}
	We have, with both $X_I$ and $X_S$ independent of $X^i$,
	\begin{align*}
	\left|\tilde{\Gamma}^N(t,X^i)-\G(t,X^i)\right|&\leq \left|\frac{1}{N} \sum_{j\in\mathfrak{I}_N}\M^N(X^i,X^j)\lambda_{-j}(t)-\overline{\lambda}^0(t)\overline{I}(0)\mathbb{E}\left(M(X^i,X_I)\Big |X^i\right)\right|\\
	&+\left|\frac{1}{N} \sum_{j\in\mathfrak{S}_N}\M^N(X^i,X^j)\lambda_{j}(t-\tilde{\tau}^N_j)-\overline{S}(0)\mathbb{E}\left(M(X^i,X_S)\overline{\lambda}(t-\tau(X_S))\big| X^i\right)\right|.
	\end{align*}
We obtain
\begin{align}
\label{im}
\mathbb{E}\left(\left|\tilde{\Gamma}^N(t,X^i)-\G(t,X^i)\right|\right)\leq \Upsilon_{1,N}^i(t)+\Upsilon_{2,N}^i(t),
\end{align}
where
	\begin{align*}
	\Upsilon_{1,N}^i(t)&=\mathbb{E}\left|\frac{1}{N} \sum_{j\in\mathfrak{I}_N}\M^N(X^i,X^j)\lambda_{-j}(t)-\overline{\lambda}^0(t)\overline{I}(0)\mathbb{E}\left(M(X^i,X_I)\Big| X^i \right)\right|,\\
	\Upsilon_{2,N}^i(t)&=\mathbb{E}\left|\frac{1}{N} \sum_{j\in\mathfrak{S}_N}\M^N(X^i,X^j)\lambda_{j}(t-\tilde\tau^N_j)-\overline{S}(0)\mathbb{E}\left(M(X^i,X_S)\overline{\lambda}(t-\tau(X_S))\Big| X^i\right)\right|\,.
	\end{align*}
	On the one hand,
	\begin{align}
	\Upsilon_{1,N}^i(t)&\leq \mathbb{E}\left(\frac{1}{N} \sum_{j\in\mathfrak{I}_N}K(X^i,X^{j}) \lambda_{-j}(t)\left|\frac{1}{\Phi\left(\int_{D}K(z,X^{j})\overline{\mu}^N(dz)\right)}-\frac{1}{\Phi\left(\int_{D}K(z,X^{j})\overline{\mu}(dz)\right)}\right|\right)\nonumber\\
	&+\mathbb{E}\left|\frac{1}{N} \sum_{j\in\mathfrak{I}_N}\left[M(X^i,X^{j})\lambda_{-j}(t)-\overline{\lambda}^0(t)\mathbb{E}\left(M(X^i,X_I)\Big| X^i\right)\right]\right|\nonumber\\
	&+\overline{\lambda}^0(t)\left|\overline{I}^N(0)-\overline{I}(0)\right|\mathbb{E}\left(M(X^i,X_I)\right)\nonumber\\
	\label{m1}
	&=\Upsilon_{1,1,N}^i(t)+\Upsilon_{1,2,N}^i(t)+\Upsilon_{1,3,N}^i(t).
	\end{align}
Exploiting in particular the inequality \eqref{LipPhi}, we have (below we assume w.l.o.g. $2\in\mathfrak{S}_N$)
	\begin{align}
\Upsilon_{1,1,N}^i(t)&\leq \lambda^*k_\Phi\mathbb{E}\left(\frac{1}{N} \sum_{j\in\mathfrak{I}_N}K(X^i,X^{j})\left|\int_{D}K(z,X^{j})\overline{\mu}^N(dz)-\int_{D}K(z,X^{j})\overline{\mu}(dz) \right|\right)\nonumber\\
&\leq C\left[ \mathbb{E} \int_D K^2(X^i,y)dy \right]^{\frac{1}{2}}\left[ \int_D  \mathbb{E} \left|\int_{D}K(z,y)\overline{\mu}^N(dz)-\int_{D}K(z,y)\overline{\mu}(dz) \right| ^2 dy\right]^{\frac{1}{2}}  \nonumber\\
&\leq  C\left[ \int_D \int_D K^2(x,y)dxdy \right]^{\frac{1}{2}}\left[ \int_D  \mathbb{E} \left|\int_{D}K(z,y)\left[\overline{\mu}^N(dz)-\overline{\mu}(dz)\right] \right| ^2 dy\right]^{\frac{1}{2}} \nonumber\\
&\leq C\left[ \int_D  \mathbb{E} \left|\int_{D}K(z,y)\left[\overline{\mu}^N(dz)-\overline{\mu}(dz)\right] \right| ^2 dy\right]^{\frac{1}{2}} \nonumber
\end{align}
It is an easy consequence of the law of large numbers that
\[	\Upsilon_{1,2,N}^i(t)\to0,\quad \text{as }N\to\infty\,.\]
\begin{align*}
\Upsilon_{1,3,N}^i(t)	&\leq c^{-\gamma}\overline{\lambda}^0(t)\left|\overline{I}^N(0)-\overline{I}(0)\right|\mathbb{E}\left(K(X^i,X_I)\right)\\
	&=c^{-\gamma}\overline{\lambda}^0(t)\left|\overline{I}^N(0)-\overline{I}(0)\right|\int_{D}\int_{D}K(x,y)\pi_S(x)\pi_I(y)dxdy\\
	&\leq C\left|\overline{I}^N(0)-\overline{I}(0)\right| 
		\end{align*}
		Summing up the above estimates, we obtain
	\begin{align}
	\label{g}
	\Upsilon_{1}^i(t)&\leq C\left[ \int_D  \mathbb{E} \left|\int_{D}K(z,y)\left[\overline{\mu}^N(dz)-\overline{\mu}(dz)\right]  \right| ^2 dy\right]^{\frac{1}{2}} + \Upsilon_{1,2,N}^i(t)+C\left|\overline{I}^N(0)-\overline{I}(0)\right|\,.
	\end{align}
	Moreover,
	\begin{align}
\Upsilon_{2,N}^i(t)&\leq \mathbb{E}\left(\frac{1}{N} \sum_{j\in\mathfrak{S}_N}K(X^i,X^{j})\lambda_{j}(t-\tilde{\tau}^N_j)\left|\frac{1}{\Phi\left(\int_{D}K(z,X^{j})\overline{\mu}^N(dz)\right)}-\frac{1}{\Phi\left(\int_{D}K(z,X^{j})\overline{\mu}(dz)\right)}\right|\right)\nonumber\\
&+\mathbb{E}\left|\frac{1}{N} \sum_{j\in\mathfrak{S}_N}M(X^i,X^{j})\lambda_{j}(t-\tilde{\tau}^N_j)-\frac{1}{N} \sum_{j\in\mathfrak{S}_N}M(X^i,X^{j})\lambda_{j}(t-\tau_j)\right|\nonumber\\
&+ \mathbb{E}\left|\frac{1}{N} \sum_{j\in\mathfrak{S}_N}\left[M(X^i,X^{j})\lambda_{j}(t-\tau_j)-\mathbb{E}\left(M(X^i,X_S)\lambda(t-\tau(X_S))\Big|X^i\right)\right]\right|\nonumber\\
&+\left|\overline{S}^N(0)-\overline{S}(0) \right|\mathbb{E}\left[M(X^i,X_S)\lambda(t-\tilde{\tau})\right]\nonumber\\
\label{m2}
&=\Upsilon_{2,1,N}^i(t)+\Upsilon_{2,2,N}^i(t)+\Upsilon_{2,3,N}^i(t)+\Upsilon_{2,4,N}^i(t)\,.
\end{align}
By the same arguments as those leading to the estimate of $\Upsilon_{1,1,N}^i$, we have
\begin{align}
\Upsilon_{2,1,N}^i(t)&\leq C\left[ \int_D  \mathbb{E} \left|\int_{D}K(z,y)\left[\overline{\mu}^N(dz)-\overline{\mu}(dz)\right] \right| ^2 dy\right]^{\frac{1}{2}} \nonumber.
\end{align}
\begin{align}
\Upsilon_{2,2,N}^i(t)&\leq\frac{C}{N}\mathbb{E}\left[\sum_{j\in\mathfrak{S}_N}K(X^i,X^{j})\left|\lambda_{j}(t-\tilde{\tau}^N_j)-\lambda_{j}(t-{\tau}_j)\right|\right]\nonumber\\
&\leq\frac{C}{N}\mathbb{E}\left[\sum_{j\in\mathfrak{S}_N}\int_D K(x,X^{j})\pi_S(x) dx \mathbf{1}_{ \left\{ \tilde{\tau}_j^N\wedge {\tau}_j\le t; \tilde{\tau}_j^N\neq {\tau}_j\right\}}\right]\nonumber\\
&\leq C\sup_{y\in D} \int_D K(x,y) dx\,\mathbb{E}\left[\frac{1}{N}\sum_{j\in\mathfrak{S}_N} \mathbf{1}_{\left\{ \tilde{\tau}_j^N\wedge {\tau}_j\le t; \tilde{\tau}_j^N\neq {\tau}_j\right\}}\right]\nonumber\\
&\leq C\mathbb{P}\left(\tilde{\tau}_2^N\wedge {\tau}_2\le t; \tilde{\tau}_2^N\neq {\tau}_2\right)\nonumber\\
&\leq C\int_{0}^{t}\mathbb{E}\left(|\tilde\Gamma^N (s,X^2)-\tilde\Gamma(s,X^2)|\right)ds\nonumber,
\end{align}
where we have used for the last inequality the facts that 
\begin{align*}
\left\{\tilde{\tau}_2^N\wedge {\tau}_2\le t; \tilde{\tau}_2^N\neq {\tau}_2\right\}&=\left\{\sup_{0 \leq t \leq T}\left|\tilde{A}^N_2(t)- A_2(t)\right|\ge1\right\},\\
\mathbb{P}\left(\sup_{0 \leq t \leq T}\left|\tilde{A}^N_2(t)- A_2(t)\right|\ge1\right)&\le
\mathbb{E}\left(\sup_{0 \leq t \leq T}\left|\tilde{A}^N_2(t)- A_2(t)\right|\right)
\end{align*}
and \eqref{t2}. Next, similarly as for $\Upsilon_{1,2,N}^i$, we deduce from the law of large numbers that
\begin{align}
\Upsilon_{2,3,N}^i(t)&\to0,\quad\text{ as }N\to\infty\nonumber.\\
\Upsilon_{2,4,N}^i(t)&\leq C\left|\overline{S}^N(0)-\overline{S}(0) \right| \nonumber
	\end{align}
	Summing up the above estimates, we obtain
	\begin{align}
	\label{g1}
	\Upsilon_{2,N}^i(t)&\leq\frac{C}{N}+C\left[ \int_D  \mathbb{E} \left|\int_{D}K(z,y)\left[\overline{\mu}^N(dz)-\overline{\mu}(dz)\right] \right| ^2 dy\right]^{\frac{1}{2}}+\Upsilon_{2,3,N}^i(t)\\
	&+C\left|\overline{S}^N(0)-\overline{S}(0)\right|+C\int_0^t \mathbb{E}\left(\left|\tilde{\Gamma}^N(s,X^2)-\G(s,X^2)\right|\right)ds\nonumber.
	\end{align}
	From \eqref{im}, \eqref{g} and \eqref{g1}, we deduce that
	\begin{align}
	\mathbb{E}\left(\left|\tilde{\Gamma}^N(t,X^i)-\G(t,X^i)\right|\right)&\leq C\left[ \int_D  \mathbb{E} \left|\int_{D}K(z,y)\left[\overline{\mu}^N(dz)-\overline{\mu}(dz)\right] \right| ^2 dy\right]^{\frac{1}{2}}\nonumber\\
	\label{t1}
	&\quad+\Upsilon_{1,2,N}^i(t)+ \Upsilon_{2,3,N}^i(t)+C\left|\overline{S}^N(0)-\overline{S}(0)\right|+C\left|\overline{I}^N(0)-\overline{I}(0)\right|\\
	&\quad+C\int_0^t \mathbb{E}\left(\left|\tilde{\Gamma}^N(s,X^2)-\G(s,X^2)\right|\right)ds\nonumber.
\end{align}
Let $\Pi^N= \displaystyle \int_D  \mathbb{E} \left|\int_{D}K(z,y)\left[\overline{\mu}^N(dz)-\overline{\mu}(dz)\right] \right| ^2 dy$, where $\overline{\mu}^N= \displaystyle\frac{1}{N} \sum_{i=1}^N \delta_{X^i},$ with $X^i$ are  i.i.d of law $\overline{\mu}$, which is the law of $X$. We have
\begin{align*}
 \left|\int_{D}K(z,y)\left[\overline{\mu}^N(dz)-\overline{\mu}(dz)\right] \right| ^2 &= \left[ \displaystyle\frac{1}{N} \sum_{i=1}^N \left[ K(X^i,y)- \mathbb{E}\left( K(X,y)\right)\right]\right]^2\\
 \mathbb{E}\left|\int_{D}K(z,y)\left[\overline{\mu}^N(dz)-\overline{\mu}(dz)\right] \right| ^2 &\leq \frac{1}{N}\mathbb{E}\left( K^2(X,y)\right)\\
  \Pi^N&\leq \frac{1}{N} \int_D\int_D K^2(x,y)dxdy.
\end{align*}

From assumption \ref{A2}, $\overline{S}^N(0)\to\overline{S}(0)$ and $\overline{I}^N(0)\to\overline{I}(0)$ a.s. as $N\to\infty$; and  $\Pi^N \to 0$ a.s. as $N\to\infty$. Hence, since moreover for any $i\in\mathfrak{S}_N$,
\[\mathbb{E}\left(\left|\tilde{\Gamma}^N(s,X^2)-\G(s,X^2)\right|\right)=\mathbb{E}\left(\left|\tilde{\Gamma}^N(s,X^i)-\G(s,X^i)\right|\right),\] it follows from
Gronwall's inequality applied to \eqref{t1} that 
\[
\mathbb{E}\left(\left|\tilde{\Gamma}^N(t,X^i)-\G(t,X^i)\right|\right)\xrightarrow [N \rightarrow \infty]{} 0.\]
The Lemma follows by combining this result with \eqref{t2}.
\end{proof} 

\subsection{Proofs of the main results}
\begin{proof}\textit{of Theorem \ref{th}.} 
	Let $\varphi \in C_b(D)$ and $T>0$ be arbitrary.
	\begin{align*}
	(\tilde{\mu}_t^{S,N},\varphi)&= (\overline{\mu}_0^{S,N},\varphi)-\frac{1}{N} \sum_{i\in\mathfrak{S}_N}\varphi(X^{i})\tilde{A}_i^N(t).
	\end{align*}
	According to the Law of Large Numbers, the sequence $(\overline{\mu}_0^{S,N},\varphi)$ converges to  $(\overline{\mu}_0^S,\varphi)\quad a.s$.\\
	 Using   Lemma \ref{Col}, as $N \rightarrow \infty$
	\begin{align*}
	\sup_{ 0\le t \leq T}\left|\frac{1}{N} \sum_{i\in\mathfrak{S}_N}\varphi(X^{i})\tilde{A}_i^N(t)-\frac{1}{N} \sum_{i\in\mathfrak{S}_N}\varphi(X^{i})A_i(t)\right|&\xrightarrow[]{Proba}0\,.
	\end{align*}
	Moreover
	\begin{align*}
	\mathbb{E}\left(\sup_{0 \leq t \leq T}\left|\varphi(X^{i})A_i(t)\right| \right) 
	\leq \|\varphi\|_{\infty}
	\end{align*}
	This last estimate allows us to use the law of large numbers in $\mathbb{D}\left(\R_+, \mathbb{R}\right)$, see Theorem 1 in \cite{rao}. 
Combining that result with Lemma \ref{Col}, we conclude that, locally uniformly in  $t$, as $N\to\infty$,
	\begin{align*}
	\frac{1}{N} \sum_{i\in\mathfrak{S}_N}\varphi(X^{i})\tilde{A}_i^N(t)&\xrightarrow[]{Proba} \overline{S}(0)\mathbb{E}\left(\int_0^t\int_0^{\infty}\mathbf{1}_{A(s^-,X_S)=0}\mathbf{1}_{u \leq \G(s,X_S)}\varphi(X_S)P(ds,du)\right)\\
	&=\overline{S}(0)\mathbb{E}\int_0^t\mathbb{P}\left(A(s^-,X_S)=0|X_S\right) \varphi(X_S)\G(s,X_S)ds\\
	&=\int_0^t\overline{S}(0)\mathbb{E}\left(\varphi (X_S)\G(s,X_S)e^{-\displaystyle \int_0^s\G(r,X_S)dr}\right)ds\\
	&=\int_0^t\int_{D}  \varphi (x)\G(s,x)\overline{\mu}_s^S(dx)ds,
	\end{align*}
	hence $(\tilde{\mu}^{S,N}_t, \varphi)\rightarrow (\overline{\mu}^{S}_t, \varphi)$ in probability locally uniformly in $t$, where
	\begin{align*}
	(\overline{\mu}_t^{S},\varphi)&= (\overline{\mu}_0^{S},\varphi)- \int_0^t\int_{D}\varphi(x)\G(s,x)\overline{\mu}_s^{S}(dx)ds.
	\end{align*}
	Next,
	\begin{align*}
	(\tilde{\mu}_t^{\mathfrak{F},N},\varphi) &= \frac{1}{N} \sum_{i\in\mathfrak{I}_N} \lambda_{-i}(t)\varphi({X^{i}}) +\frac{1}{N} \sum_{i\in\mathfrak{S}_N} \lambda_{i}(t- \tilde{\tau}_i^N)\varphi({X^{i}}).
	\end{align*}
	On the one hand,  as $N\to\infty$,
	\begin{align*}
	\frac{1}{N} \sum_{i\in\mathfrak{I}_N} \lambda_{-i}(t)\varphi({X^{-i}})&\xrightarrow{a.s}\overline{I}(0)\mathbb{E}\left(\lambda_{-1}(t)\varphi({X_I})\right)
	=\overline{\lambda}^0(t)(\overline{\mu}^I_0,\varphi),\\
		 \text{where}\quad \overline{\mu}^I_0(dx)=\overline{I}(0) \pi_I(x)dx\,.
	\end{align*} 
Moreover, $\displaystyle \frac{1}{N} \sum_{i\in\mathfrak{S}_N} \lambda_{i}(t- \tilde{\tau}_i^N)\varphi({X^{i}})=
\frac{1}{N} \sum_{i\in\mathfrak{S}_N} \int_0^t \lambda_{i}(t- s)\varphi({X^{i}})d\tilde{A}_i^N(s)$.\\ 
According to a variant of Lemma \ref{Col}, as $N\rightarrow \infty$\\
	\[\sup_{ t \leq T}\left| \frac{1}{N} \sum_{i\in\mathfrak{S}_N} \int_0^t \lambda_{i}(t- s)\varphi({X^{i}})d\tilde{A}_i^N(s)-\frac{1}{N} \sum_{i\in\mathfrak{S}_N} \int_0^t \lambda_{i}(t- s)\varphi({X^{i}})dA_i(s)\right|\xrightarrow{Proba}0.\] 
	Indeed,
	\[\sup_{ t \leq T}\left| \frac{1}{N} \sum_{i\in\mathfrak{S}_N} \int_0^t \lambda_{i}(t- s)\varphi({X^{i}})\left(d\tilde{A}_i^N(s)-dA_i(s)\right)\right|\leq \frac{\lambda^*\|\varphi\|_{\infty}}{N}\sum_{i\in\mathfrak{S}_N}\sup_{0 \leq t \leq T}\left| \tilde{A}^N_i(t)-A_i(t)\right|.\]
	In addition,
	\begin{align*}
	\mathbb{E}\left(\sup_{0 \leq t \leq T}\left|\int_0^t  \lambda_{i}(t- s)\varphi({X^{i}})dA_i(s)\right|\right)
	&\leq \lambda^* \|\varphi\|_{\infty}.
	\end{align*}
	Hence combining the above result with the law of large numbers in $\mathbb{D}\left(\R_+, \mathbb{R}\right)$, we deduce that, in probability, locally uniformly in $t$, as $N\to\infty$,
	\begin{align*}
	\frac{1}{N} \sum_{i\in\mathfrak{S}_N}&\int_0^t\lambda_{i}(t- s)\varphi({X^{i}})dA_i(s)\\
	&\xrightarrow{\mathbb{P}} \overline{S}(0)\mathbb{E}\left(\int_0^t\lambda(t- s)\int_0^{\infty}\mathbf{1}_{A(s^-,X_S)=0}\mathbf{1}_{u \leq \G(s,X_S)}\varphi(X_S)P(ds,du)\right)\\
	&=\overline{S}(0)\int_0^t\overline{\lambda}(t-s)\mathbb{E}\left[\mathbb{P}\left(A(s^-,X_S)=0|X_S\right) \varphi(X_S)\G(s,X_S)\right]ds\\
	&=\int_0^t\overline{\lambda}(t-s)\overline{S}(0)\mathbb{E}\left(\varphi (X_S)\G(s,X_S)e^{-\displaystyle \int_0^s\G(r,X_S)dr}\right)ds\\
	&=\int_0^t\overline{\lambda}(t-s)\int_{D}  \varphi (x)\G(s,x)\overline{\mu}_s^S(dx)ds\,.
	\end{align*}
We thus obtain that $(\tilde{\mu}_t^{\mathfrak{F},N},\varphi) \rightarrow (\overline{\mu}_t^{\mathfrak{F}},\varphi)$ in probability locally uniformly in t, where
	\begin{align*}
	(\overline{\mu}_t^{\mathfrak{F}},\varphi) &=\overline{\lambda}^0(t)(\overline{\mu}^I_0,\varphi) +\int_0^t\overline{\lambda}(t- s)\int_{D} \varphi(x)\G(s,x)\overline{\mu}_s^{S}(dx)ds. 
	\end{align*}
	Now,
	\begin{align*}
	(\tilde{\mu}_t^{I,N},\varphi)&= (\overline{\mu}_0^{I,N},\varphi)+\frac{1}{N} \sum_{i\in\mathfrak{S}_N}\varphi(X^{i})\tilde{A}^N_i(t)\\&\quad- \frac{1}{N}\sum_{i\in\mathfrak{I}_N}\varphi(X^{i})\mathbf{1}_{\eta_{-i}\leq t}- \frac{1}{N} \sum_{i\in\mathfrak{S}_N}\varphi(X^{i})\int_0^t\mathbf{1}_{\eta_i\leq t-s}d\tilde{A}^N_i(s)
	\end{align*}
	According to the Law of Large Numbers:
	\begin{align*}
	(\overline{\mu}_0^{I,N},\varphi)&\xrightarrow[]{a.s.}(\overline{\mu}_0^I,\varphi),\\
	\frac{1}{N}\sum_{i\in\mathfrak{I}_N}\varphi(X^{i})\mathbf{1}_{\eta_{-i}\leq t}&\xrightarrow[]{a.s.} \overline I(0)\mathbb{E}\left(\varphi(X_I)\mathbf{1}_{\eta_{-1}\leq t}\right)
	=F_0(t)(\overline{\mu}_0^I,\varphi)\,.
	\end{align*}
	 
	From the previous results,
	\begin{align*}
	\frac{1}{N} \sum_{i\in\mathfrak{S}_N}\varphi(X^{i})\tilde{A}^N_i(t)\xrightarrow[]{\mathbb{P}}\int_0^t\int_{D}\varphi(x)\G(s,x)\overline{\mu}_s^{S}(dx)ds
	\end{align*}
	According to a variant of Lemma \ref{Col}, as $N\rightarrow \infty$,
	\[\sup_{ t \leq T}\left| \frac{1}{N} \sum_{i\in\mathfrak{S}_N}\varphi(X^{i})\int_0^t\mathbf{1}_{\eta_i\leq t-s}d\tilde{A}^N_i(s)-\frac{1}{N} \sum_{i\in\mathfrak{S}_N}\varphi(X^{i})\int_0^t\mathbf{1}_{\eta_i\leq t-s}dA_i(s)\right|\xrightarrow{Proba}0.\]
	 Indeed, 
	\begin{align*}
	\sup_{ t \leq T}&\left| \frac{1}{N} \sum_{i\in\mathfrak{S}_N}\varphi(X^{i})\int_0^t\mathbf{1}_{\eta_i\leq t-s}d\tilde{A}^N_i(s)-\frac{1}{N} \sum_{i\in\mathfrak{S}_N}\varphi(X^{i})\int_0^t\mathbf{1}_{\eta_i\leq t-s}dA_i(s)\right|
	\\ &\leq \frac{\|\varphi\|_{\infty}}{N}\sum_{i\in\mathfrak{S}_N}\left(\sup_{0 \leq t \leq T}\left| \tilde{A}^N_i(t)-A_i(t)\right|\right)\,.
	\end{align*}
	In addition
	\begin{align*}
	\mathbb{E}\left(\sup_{0 \leq t \leq T}\left|\int_0^t\mathbf{1}_{\eta_i\leq t-s}\varphi({X^{i}})dA_i(s)\right|\right)
	\leq \|\varphi\|_{\infty}\,.
	\end{align*}
		Applying the Law of Large Numbers in $\mathbb{D}\left(\R_+, \mathbb{R}\right)$, we deduce that, locally uniformly in $t$,
		as $N\to\infty$,
	\begin{align*}
	\frac{1}{N} \sum_{i\in\mathfrak{S}_N}&\int_0^t\mathbf{1}_{\eta_i\leq t-s}\varphi({X^{i}})dA_i(s)\\
	 &\xrightarrow{Proba}\overline{S}(0)\mathbb{E}\left(\int_0^t\mathbf{1}_{\eta_1\leq t-s}\int_0^{\infty}\mathbf{1}_{A(s^-,X_S)=0}\mathbf{1}_{u \leq \G(s,X_S)}\varphi(X_S)P(ds,du)\right)\\
	&=\overline{S}(0)\int_0^tF(t-s)\mathbb{E}\left[\mathbb{P}\left(A(s,X_S)=0|X_S\right) \varphi(X_S)\G(s,X_S)\right]ds\\
	&=\int_0^tF(t-s)\overline{S}(0)\mathbb{E}\left(\varphi (X_S)\G(s,X_S)e^{-\displaystyle \int_0^s\G(r,X_S)dr}\right)ds\\
	&=\int_0^tF(t-s)\int_{D}  \varphi (x)\G(s,x)\overline{\mu}_s^S(dx)ds\,.
	\end{align*}
Combining the above results, we deduce that $(\tilde{\mu}_t^{N,I},\varphi)\rightarrow (\overline{\mu}_t^{I},\varphi)$ in probability locally uniformly in $t$, where
	\begin{align*}
	(\overline{\mu}_t^{I},\varphi)&=(\overline{\mu}_0^{I},\varphi)+\int_0^t\int_{D} \varphi(x)\G(s,x)\overline{\mu}_s^{S}(dx)ds\\
	&\quad-(\overline{\mu}_0^{I},\varphi) F_0(t)-\int_0^t F(t-s)\int_{D} \varphi (x)\G(s,x)\tilde{\mu}_s^S(dx)ds,\nonumber
	\end{align*}
	which can be rewritten as:
	\begin{align*}
	(\overline{\mu}_t^{I},\varphi)&= (\overline{\mu}_0^{I},\varphi) F^c_0(t)+\int_0^tF^c(t-s)\int_{D} \varphi (x)\G(s,x)\overline{\mu}_s^S(dx)ds\,.
	\end{align*}
	
	We argue analogously about  $\tilde{\mu}^{R,N}$ and obtain that 
	$(\tilde{\mu}_t^{R,N},\varphi)\rightarrow(\overline{\mu}_t^{R},\varphi)$ in probability locally uniformly in $t$, where
	\begin{align*}
	(\overline{\mu}_t^{R},\varphi)&=(\overline{\mu}_0^{R},\varphi)+(\overline{\mu}_0^{I},\varphi) F_0(t)+\int_0^tF(t-s)\int_{D}\varphi (x)\G(s,x)\overline{\mu}_s^S(dx)ds.
	\end{align*}
\end{proof}

Let us now establish a technical result, which will be useful for the last step of our proof.

For the next Lemma, we shall need the following construction.
Consider the partition of $\R^d$ made of translates of the hypercube $(0,\frac{a}{\sqrt{d}}]^d$, where the value of $a$ will be specified in the proof of the next Lemma.
Let $(\mathfrak{A}_1, \mathfrak{A}_2,\ldots, \mathfrak{A}_k)$ denote the subset of the elements of that partition which are contained in $D$. 

In the next statement, $r$ is the constant which appears in assumptions \ref{A1} and \ref{A2}.
\begin{lemma}
\label{impor}
For all $y\in D$, there exists $1\leq j\leq k$ such that $\mathfrak{A}_j \subset  B(y,r)$.
\end{lemma}
\begin{proof}
Let $y\in D$. According to the hypothesis, $\exists \alpha>0$ such that $C(y,\ell_y, \alpha)\cap B(y,r) \subset D$.\\
 Let  $u(y):=y+\frac{r}{1+\sin(\alpha)}\ell_y$ and $a:=\frac{r\sin(\alpha)}{1+\sin(\alpha)}.$ Then  $B(u(y),a)\subset C(y,\ell_y, \alpha)\cap B(y,r)  $. In fact,  $\forall z\in B(u(y),a)$ we have
 \[ \|y-z\|\leq \|y-u(y)\|+\|u(y)-z\| \leq \frac{r}{1+\sin(\alpha)}+\frac{r\sin(\alpha)}{1+\sin(\alpha)}=r,\ \text{ hence }z \in B(y,r).\]
The minimum distance from $u(y)$ to the boundary of the cone $C(y,\ell_y,\alpha)$ is
\[ \|u(y)-y\|\sin (\alpha)= \frac{r\sin (\alpha)}{1+\sin(\alpha)}\geq \|z-u(y)\|,\text{ hence }z \in C(y,\ell_y, \alpha).\]
Let $(e_1,\cdots,e_d),$ denote an orthonormal basis of $\R^d$. The hypercube $H(y):=\left\{ u(y)+\displaystyle \sum_{i=1}^d t_i e_i, |t_i|\leq \frac{a}{\sqrt{d}}\right\}$ is contained in $B(u(y),a)\subset D$. We also have $H(y)\subset B(y,r)$. It is plain that there exists $1\le j\le k$ such that $\mathfrak{A}_j\subset H(y)$: we can choose the unique $j$ s.t. $u(y)\in\mathfrak{A}_j$.
\end{proof}
\begin{proof}\textit{of Theorem \ref{theo1}}.
We have,
\begin{align*}
\inf_{y\in D}\int_{D} K(z,y)\overline{\mu}(dz)&\geq \underline{c} \inf_{y\in D}\int_{C(y,\ell_y,\alpha)\cap B(y,r)} \overline{\mu}(dz)\\
                     &\geq  \underline{c} \inf_{1\leq i\leq k} \overline{\mu}(\mathfrak{A}_i)\\
                     &\geq  \underline{c} a^d d^{-d/2} \inf_{z\in D}\overline{\mu}(z):= \theta
\end{align*}
We now choose $\delta=\theta/2$ and define $\displaystyle \Omega_N:= \left\{ \omega, \inf_{y\in D}\int_{D} K(z,y)\overline{\mu}^N(dz)> \delta\right\}$. \\
 We remark that on $\Omega_N$, $\left(\tilde{\mu}^{S,N},\tilde{\mu}^{\mathfrak{F},N},\tilde{\mu}^{I,N}, \tilde{\mu}^{R,N}\right)=\left(\overline{\mu}^{S,N},\overline{\mu}^{\mathfrak{F},N},\overline{\mu}^{I,N}, \overline{\mu}^{R,N}\right)$. Hence  clearly Theorem \ref{theo1}  will follow from Theorem \ref{th} if we prove that $\P(\Omega_N)\to1$, as $N\to\infty$. Fix $y\in D$ and
let $j$ be such that $\mathfrak{A}_j \subset  B(y,r)$. On the event that all the $X^i$, $ i\ge 1$, are distinct, which happens a.s., for each $N\ge1$,
\begin{align*}
\int_{D} K(z,y)\overline {\mu}^N(dz)&\geq \underline{c}\, \overline{\mu}^N (\mathfrak{A}_i\backslash\{y\})\\
&\ge\underline{c}\,\left(\overline{\mu}^N (\mathfrak{A}_j)-\frac{1}{N}\right)\, ,\ \text{then}\\
\inf_{y\in D}\int_{D} K(z,y)\overline {\mu}^N(dz)&\ge\underline{c}\,\left(\inf_{1\le j\le k}\overline{\mu}^N (\mathfrak{A}_j)-\frac{1}{N}\right)\,.
\end{align*}
Consequently
\[\liminf_{N\rightarrow \infty}\inf_{y\in D}\int_{D} K(z,y)\overline {\mu}^N(dz)\geq \underline{c} \inf_{1\leq i\leq k} \overline{\mu} (\mathfrak{A}_i)\ge\theta=2\delta.
\]
Hence $\P(\Omega_N)\to1$, as $N\to\infty$.
\end{proof}

		\paragraph{\bf Acknowledgement} The authors wish to thank an anonymous Referee, whose detailed comments helped the authors to improve significantly the original version of this paper, in particular concerning the proof of Lemma \ref{lll}.

\end{document}